\documentclass[10pt]{amsart}
\usepackage{amsmath,amssymb,amsthm, epsfig}
\usepackage{hyperref}
\usepackage{amsmath}
\usepackage{amssymb}
\usepackage{hyperref}
\usepackage{ulem}

\title[Optimization problem for degenerate quasilinear operator]{An optimization problem with free boundary governed by
a degenerate quasilinear operator.}
\author{Krerley Oliveira}
\address{Universidade Federal de Alagoas, Departamento de Matem\'{a}tica, 57072-900 Macei\'{o}-Alagoas}
\email{krerley@mat.ufal.br}
\author{Eduardo V. Teixeira}
\address{Department of Mathematics, Rutgers University, Hill Center-Busch Campus, 110 Frelinghuysen Road, Piscataway, NJ 08854-8019}
\email{teixeira@math.rutgers.edu}

\date{}

\def \R {\mathbb{R}}

\def \div {\textrm{div}}

\newtheorem{theorem}{Theorem}[section]
\newtheorem{lemma}[theorem]{Lemma}

\theoremstyle{definition}

\theoremstyle{remark}

\numberwithin{equation}{section}

\begin{document}
\maketitle

\begin{abstract}
In this paper we study the existence, regularity and geometric properties of
an optimal configuration to a free boundary optimization problem
governed by the $p$-Laplacian.

\end{abstract}

\section{Introduction}

Let $D$ be a bounded domain in  $\mathbb{R}^n$ and  $\varphi$ a positive
function defined on it (the temperature distribution of the
body $D$). A classical minimization problem in heat conduction
asks for the best way of insulating the body $D$, with a prescribed amount of insulating material in a stationary situation.
This situation also models problems in electrostatic, potential
flow in fluid mechanics among others.  \par

The mathematical description of this problem is
as follows: fixed a number, $\gamma > 0$ (quantity of insulation material), for each domain $\Omega$ surrounding $D$, such that $|\Omega \setminus D| = \gamma$, we consider
the potential $u$ associated to the configuration $\Omega$, i.e.,
the harmonic function in $\Omega \setminus D$, taking boundary data equal to
$\varphi$ on  $\partial D$  and $0$ on $\partial\Omega$.
The flow of heat (quantity to be minimized) corresponding to the configuration $\Omega$, is
given by a nonlocal monotone operator given by,
$$
J(\Omega) := \int_{\partial D} \Gamma\big (x, u_\mu(x) \big ) dS,
$$
where $\mu$ is the inward normal vector defined on $\partial D$.
The function $\Gamma \colon \partial D \times \mathbb{R} \to
\mathbb{R}$ is assumed to be convex and increasing, on $u_\mu$ and
continuous on $x$. Important examples are, $\Gamma(t) = t$
(classical heat conduction problem), $\Gamma(t) = t^p$ (optimal
configurations in electrostatics), $\Gamma(x,t) = \max\{t, C(x)\}$
(problems is the material sciences).  \par In \cite{T}, the second
author studied the problem of minimizing $J(\Omega)$ among all
configurations $\Omega$ such that, say, $|\Omega \setminus D| =
1$, where $|A|$ is the volume of the set $A$. This optimization
problem with linear heat flux, i.e., $\Gamma(x,t) = t$, was
studied in \cite{AAC} and \cite{ACS}. Qualitative geometric
properties of the free boundary, namely symmetry, uniqueness and
full regularity of the free boundary, was explored in \cite{T1}.
\par In this present paper, we turn our attention to this problem
when we allow the temperature itself deform the medium. We assume
the influence of the temperature distribution on the medium is
proportional to the magnitude of its gradient. These
considerations lead us to study this optimization problem when
temperature distribution is governed by the $p$-Laplacian. In
other words, the variational problem we are interested in is, for
$1<p<\infty$
\begin{equation} \label{Q}
\textrm{Minimize } \left \{
    \begin{array}{c}
        J(u) := \displaystyle\int_{\partial D} \Gamma(x, u_\mu)
 dS ~\Big | ~ u \colon D^C \to \R, ~ u = \varphi
        \textrm{ on } \partial D, \\
         ~\Delta_p u = 0  \textrm{ in }\{u>0\} \textrm{ and } |\{u>0\}|
        = 1
        \end{array} \right \}
\end{equation}

Here, $\Delta_p u := \div \big ( |\nabla u |^{p-2}\nabla u \big
)$. Several new difficulties appear when dealing with the
nonlinear operator $\Delta_p$. One of its main difficult lies in
the fact that the $p$-Laplacian is not uniformly elliptic.  \par

In analogy with the linear case for the Laplacian operator, in
this paper we shall restrict ourselves to the heat flux given by $\Gamma(x,t) = t^{p-1}$.
In our physical considerations, we will assume the body to be insulated has much smaller volume than the quantity of insulation material. This leads us to consider a constant temperature distribution, say, $\varphi \equiv 1$. With these assumptions, problem
(\ref{Q}) can be reformulated in terms of the following equivalent
version of it:
\begin{equation} \label{P}
\textrm{Minimize } \left \{
    \begin{array}{c}
        J(u) := \displaystyle\int_{D^C} |\nabla u|^p
 dx ~\Big | ~ u \colon D^C \to \R, ~ u = \varphi
        \textrm{ on } \partial D, \\
         ~\Delta_p u = 0  \textrm{ in }\{u>0\} \textrm{ and } |\{u>0\}|
        = 1
        \end{array} \right \}.
\end{equation}


 In problem (\ref{P}), we may assume $\varphi$ to be a mere continuous and
positive function. When $\varphi \equiv 1$ in problem (\ref{P}), such a problem is equivalent to problem
({\ref{Q}) with $\Gamma(x,t) = t^{p-1}$. In this paper we will only deal with problem (\ref{P}). We hope to turn
our attention to problem (\ref{Q}) in its full generality in future research.\par

From the mathematical point of view, our approach is motivated by
recent advances on the free boundary regularity theory for minimum
problems with a variable domain of integration involving degenerate
quasilinear operators. Namely, D. Danielli and A. Petrosyan, in
\cite{DP}, have recently extended the celebrated work of H. Alt and L. Caffarelli
\cite{AC}, for the $p$-Laplacian operator. Our regularity results
will relay on suitable modifications of the arguments in \cite{DP}.
Furthermore, we shall establish a free boundary condition that will
relate our optimization problem with Bernoulli type problems,
similar to the ones studied in \cite{HS}, \cite{DPS} and \cite{DP}.

\section{Mathematical fundaments of the physical problem}

In this section we shall introduce the main mathematical tools we
shall use throughout the whole paper. Throughout the article, $1 <
p < \infty$ and $\Delta_p u $ stands for the $p$-Laplacian
operator
$$
    \Delta_p u := \div \big (|\nabla|^{p-2}\nabla u \big ).
$$
Let $U$ be a domain in $\mathbb{R}^N$. Let us recall that for any
$\xi \in W^{1,p}(U)$, $\Delta_p \xi \in \big [W_0^{1,p}(U) \big
]^{*}$ and
$$
    \langle \Delta_p \xi, \upsilon \rangle = \int_U |\nabla
    \xi|^{p-2} \nabla \xi \cdot \nabla \upsilon dx.
$$
Problem (\ref{P}) presents several difficulties from the
mathematical point of view. Our strategy will be to study a penalized version
of it, which is sort of a weak formulation of the problem. The
idea is to grapple with the difficulty of volume constraint, which
is very unstable under limits and makes perturbation arguments quite hard. \par

From now on, we denote by $V$ the following set
$$
    V := \left \{ u \in W^{1,p}(D^C) ~\Big | ~ u = \varphi
    \textrm{ on } \partial D \right \}.
$$
The penalized problem is stated as follows: Let $ \varepsilon > 0$
be fixed. We consider the function
$$
f_\varepsilon := \left \{
\begin{array}{rll}
    1 + \dfrac{1}{\varepsilon} (t-1) & \textrm{ if }  t \ge 1& \\
    1 + \varepsilon (t-1) & \textrm{ otherwise.} &
\end{array} \right.
$$
We then define the penalized functional as
\begin{equation}\label{Je}
J_\varepsilon(u) := \int_{D^C} |\nabla u|^{p} dx +
f_\varepsilon\Big ( | \{u>0\} | \Big ).
\end{equation}
For the moment, we shall be interested in the following minimization
problem
\begin{equation} \label{Pe}
\min\limits_{u \in V} J_\varepsilon(u).
\end{equation}

For latter use, given a $\gamma$-bilipschitz function $f$ we
consider the similar penalized problem of (\ref{Pe}) given by the
functional:
 \begin{equation}\label{Pf}
 J_f:= \int_{D^C} |\nabla u|^{p} dx +
f\Big ( | \{u>0\} | \Big ).
\end{equation}

\section{Properties of solutions of Problem (\ref{Pe})}

In this section we shall derive existence of a minimizer for the penalized problem as well as some important nondegeneracy conditions, such as optimal regularity and linear growth away from the free boundary. The proof of Theorem \ref{Existence} will be developed throughout this section. At the end, we shall be able to state a representation theorem that will be crucial to study further regularity properties of the free boundary.
\begin{theorem}\label{Existence} For each $\varepsilon > 0$ fixed, there
exists a minimizer $u_\varepsilon\in V$ for the functional
$J_\varepsilon$. Furthermore
\begin{enumerate}
\item $u_\varepsilon \ge 0$.
\item $\Delta_p u_\varepsilon$ is a nonnegative Radon measure
supported on $\partial \{ u_\varepsilon > 0 \}$. In particular,
$$
    \Delta_p u_\varepsilon = 0  \textrm{ in } \{ x \in D^C ~\Big | ~
u_\varepsilon (x) > 0 \}.
$$
\item Any minimizer $u_\varepsilon$ of problem (\ref{Pe}) is  Lipschitz continuous and for
any compact $\mathcal{K} \subset D^C$ there exists a constant $K=K(\mathcal{K}, p,n,\varepsilon,D)$ such that
$\| u_\varepsilon\|_{\textrm{Lip}(\mathcal{K})} \le K$.
\item The function $u_\varepsilon$ growths linearly way from the
free boundary, i.e., for any compact $\mathcal{K} \subset D^C$,
there exist positive constants $c, ~C$, depending on dimension,
$\mathcal{K}$, $p$, $D$ and $\varepsilon$, such that
$$
    c \textrm{dist}\big ( x, \partial \{ u_\varepsilon > 0 \} \big
    ) \le u_\varepsilon(x) \le C \textrm{dist}\big ( x, \partial \{ u_\varepsilon > 0 \} \big
    ), ~ \forall x \in \mathcal{K}.
$$
\item The free boundary is uniformly dense, i.e., for any
compact $\mathcal{K} \subset D^C$ fixed, there exist a constant $c=
c(\mathcal{K}, p, n,\varepsilon, D)$ with $0< c < \frac{1}{2}$, such
that
$$
    c < \dfrac{| B_r \cap \{ u > 0 \}|}{|B_r|} \le 1-c,
$$
for any ball  $B_r = B_r(x)$ centered at some point $x \in
\partial \{u_\varepsilon > 0 \} \cap \mathcal{K}$.
\end{enumerate}
\end{theorem}
\begin{proof}
The existence of a minimizer $u_\varepsilon$ for problem (\ref{Pe}) follows easily from the fact that for any
minimizing sequence $u^k_\varepsilon$, we may assume

\begin{itemize}
\item $\nabla u^k_\varepsilon \rightharpoonup \nabla u_\varepsilon$
in $L^p(D^C)$
\item $u^k_\varepsilon \to u_\varepsilon$ almost
everywhere in $D^C$,
\end{itemize}
for some $u_\varepsilon \in V$. Thus,
$$
    \begin{array}{c}
       \displaystyle \int_{D^C} |\nabla u_\varepsilon|^p dx \le \liminf\limits_{k
        \to \infty}  \displaystyle \int_{D^C} |\nabla u^k_\varepsilon|^p dx \\
        \textrm{and} \\
        |\{ u_\varepsilon > 0 \}| \le \liminf\limits_{k
        \to \infty} |\{ u^k_\varepsilon > 0 \}|.
    \end{array}
$$
Since $f_\varepsilon$ is a continuous and increasing function, we
obtain
$$
    J_\varepsilon(u_\varepsilon) \le \liminf\limits_{k
    \to \infty} J(u_\varepsilon^k).
$$
Clearly, $u_\varepsilon \ge 0$, otherwise, $J_\varepsilon\big (
~(u_\varepsilon)^{+} \big ) < J_\varepsilon\big ( u_\varepsilon\big
)$. \par

Observe that if $u_\varepsilon$ is a minimizer, then $J(u_\varepsilon) \leq J(u_\varepsilon-\epsilon \eta)$, for
every $\epsilon
>0$ and nonnegative $\eta \in C^\infty_0(D^c)$. Since $\{u_\varepsilon>\epsilon
\eta\} \subset \{u\varepsilon>0\}$ and $f$ is increasing, we have
that $f \big (\{u_\varepsilon>\epsilon \eta\} \big ) \leq  f\big
(\{u_\varepsilon>0\}\big)$  and consequently, we have that $\Delta_p
u_\varepsilon$ is a nonnegative Radon measure supported on
$\partial \{u_\varepsilon >0\}.$

Now, we explain the main track and the necessary changes in order to
obtain items $(3)$, $(4)$ and $(5)$ . We will follow the lines of
\cite{DP} establishing a sequence of lemmas that are analogous of
those in this refereed work. Instead of proving these lemmas with all the details, we
shall restrict ourselves to enunciate them and sketch their proofs by pointing
out the necessary modifications and sending to \cite{DP} for further details.

\textit{Up to the end of this section, we  fix some
$\gamma$-bilipschitz function $f$.} Consider $u$ a minimizer of the
problem \ref{Pf} in some ball $B$. Denote by $v$ the solution of the
Dirichlet problem
\begin{equation} \label{Sol. p-Laplacian}
    \left \{
        \begin{array}{rlcl}
            \Delta_p v&=& 0 &\text{ in }B\\
            v& = & u &\text{ on }\partial B.
        \end{array}
    \right.
\end{equation}
Following the beginning of Section $3$ in \cite{DP} we obtain:

\begin{lemma} \label{l.contagradiente}
There exists a constant $C=C(n,p,\gamma)$ such that

\begin{equation}
\begin{split}
&\int_B |\nabla (u-v)|^p \leq  C |\{u>0\}\cap B| \text{ for }p \geq 2 \text{ and }\\
 \int_B |\nabla (u&-v)|^p \leq  C |\{u>0\}\cap B|^{\frac{p}{2}} \big ( \int_B |\nabla u|^p \big )^{1-\frac{p}{2}} \text{ for } 1<p \leq 2.
\end{split}
\end{equation}
Moreover, the constant $C$ goes to zero when $\gamma$ goes to zero.
\end{lemma}

The next Lemma is the analogous of Lemma~3.1 in \cite{DP}, with a similar proof.
\begin{lemma}\label{l.holder}
Let $f$ be a  given $\gamma$- Lipschitz function and $u$ be a bounded minimizer of Problem \ref{Pf} in $B_1$. Then, $u \in C^{\alpha}$ in $B_{\frac{7}{8}}$ for some $\alpha=\alpha(n,p) \in (0,1)$ and
 $$
  \|u\|_{C^{\alpha}(B_{\frac{7}{8}})} \leq C(n,p,\|u\|_{L^\infty(B_1)}, \gamma).
  $$
\end{lemma}
Following Lemma~3.2 in \cite{DP}, we obtain
\begin{lemma}\label{l.limitado}
Let $u$ be a bounded  local minimizer of Problem \ref{Pf} in $B_1$ with $u(0)=0$. Then, there
exists a constant $C=C(n,p,\gamma)$ such that
$$
\|u\|_{L^\infty(B_{ \frac{1}{4}})} \leq C.
$$
\end{lemma}
\begin{proof} Assume, by contradiction, there exists a sequence $u_k$ of bounded local minimizers of Problem~\ref{Pf}
in $B_1$ with $u_k(0)=0$ and $\max_{B_{\frac{1}{4}}} u_k>k.$  In this case, we define:
$$
    d_k(x)=d(x, \partial\{u_k>0\})
$$
and
$$
    \mathcal{O}_k = \left \{ x\in B_1; d_k(x)\leq \frac{1-\|x\|}{3} \right \}.
$$
Observe that, since $u_k(0)=0$, $d_k(x) < \|x\|$ for every $x \in B_1$. On the other hand,  $(1-\|x\|)/3 \geq
1/4$ for every $x \in B_{1/4}$. From this inequality we conclude that $B_{1/4} \subset \mathcal{O}_k$. Now, define
$$
m_k := \max_{\mathcal{O}_k} (1-\|x\|)u_k(x) \geq \frac{3}{4} \max_{B_{1/4}} u_k(x) > \frac{3}{4} k.
$$
Consider any maximum point $x_k \in \mathcal{O}_k$ of the function $(1-\|x\|)u_k(x)$ and observe that
\begin{equation}\label{e.1}
u_k(x_k)= \frac{m_k}{1-\|x\|} > \frac{3}{4}k.
\end{equation}
Denote by $y_k$  any point in $\partial\{u_k>0\}$
such that $d_k(x_k)=\|y_k-x_k\|$ and define $\delta_k := \|x_k-y_k\|$. Since $ x_k \in \mathcal{O}_k$ we have
that $\delta_k \leq \frac{(1-\|x\|)}{3}$. Then, for every $z \in B_{2\delta_k}(y_k)$:
$$
\|z\| \leq \|y_k\| + 2 \delta_k \leq \|x_k\| + 3 \delta_k \leq \|x_k\| +  (1-\|x_k\|) \leq 1.
$$
From this inequality, we have that $B_{2\delta_k}(y_k) \subset B_1$. Now, we claim that $B_{\delta_k/2}(y_k)
\subset \mathcal{O}_k.$ In fact, if $z \in B_{\delta_k/2}(y_k)$ we have that:
$$
d_k(z) \leq \frac{\delta_k}{2} \leq \frac{1-\|y_k\| + \frac{\delta_k}{2}}{3} \leq \frac{(1-\|z\|)}{3},
$$ where the second inequality is a consequence of $\|y_k\|+\delta_k \leq 1$. Moreover, for $z\in B_{\delta/2}(y_k)$
$$
(1-\|z\|) \geq (1-\|x_k\|)-\|x_k-z\| \geq  (1-\|x_k\|)-\frac{3}{2}\delta_k \geq \frac{(1-\|x_k\|)}{2}.
$$
As a consequence of this inequality,
\begin{equation}\label{e.2}
\max\limits_{\overline{B}_{\delta/2}(y_k)} u_k \leq 2u_k(x_k).
\end{equation}

By the definition of $\delta_k$, we have that $B_{\delta_k}(x_k) \subset \{u_k>0\}.$ Recall from item
(2) of Theorem~\ref{Existence} that $\Delta_p u_k = 0$ in $B_{\delta_k}(x_k)$. By the Harnack
Inequality for $p$-harmonic functions, we may conclude that there exists a constant $c=c(n,p)>0$ such that
\begin{equation}\label{e.h}
\min\limits_{\overline{B}_{3\delta/4}(x_k)} u_k \geq cu_k(x_k).
\end{equation}
In particular,
\begin{equation}\label{e.3}
\max\limits_{\overline{B}_{\delta/4}(y_k)} u_k \geq cu_k(x_k).
\end{equation}
Consider the following scaling of $u_k$
\begin{equation}\label{e.4}
 w_k(x)=\frac{u_k(y_k+\frac{\delta}{2} x)}{u_k(x_k)}, \text{ for } x \in B_1
\end{equation}

Observe that since $u$ is a local minimizer of Problem~\ref{Pf}, $w_k$ is a minimizer of the analogous
problem replacing $f_\epsilon$ by $f_k = f_\epsilon/u_k(x_k)$. In other words, $w_k$ is a local minimizer of
$$
J_k(w) = \int_{B_{3/4}} \|\nabla w\|^p + f_k(|\{w>0\}|).
$$
Now, we denote by $v_k$ the solution of $\Delta_p v_k = 0$ in $B_{3/4}$ and $v_k-w_k \in W_0^{1,p}(B_{3/4})$.
Initially we observe that $\|f_k\|_{Lip}$  converges uniformly to zero, as $k \rightarrow \infty$. Since $w_k$ is a
minimizer of $J_k$, by Lemma~\ref{l.contagradiente} we guarantee the existence of a sequence of positive numbers $C_k$, that converges to zero
as $k \to \infty$, such that
\begin{equation}\label{e.5}
 \int_{B_{3/4}} \|\nabla (w_k- v_k)\|^p  \leq   C_k.
\end{equation}
By estimates~(\ref{e.2}) and (\ref{e.3}) we have that
\begin{equation*}
\max\limits_{\overline{B}_1} w_k \leq 2\text{, } \max\limits_{\overline{B}_{1/2}} w_k \geq c\text{ and }
w_k(0)=0.
\end{equation*}
Now notice that $w_k$ is uniformly bounded, thus from Lemma~\ref{l.holder} we conclude $w_k$ and $v_k$ are
uniformly $C^\alpha$ in $B_{5/8}$. By compactness, we may assume (passing to a subsequence, if necessary) $w_k \rightarrow w_0$ and $v_k \rightarrow v_0$
uniformly on $B_{5/8}$. Equation~(\ref{e.5}) implies that $w_0=v_0+K$ in $B_{5/8}$. Since $v_k \rightarrow v_0$, we have that $\Delta_p v_0 = 0$ and this implies that $\Delta_p w_0 =
0$.  By the strong minimum principle, we have that $w_0=0$ in $B_{5/8}$, because $w_0\geq 0$ and $w_0(0)=0$. It is
a contradiction with $\max\limits_{B_{1/2}} w_0>c>0$. This finishes the proof.
\end{proof}
Now, we show how one can derive Lipschitz continuity of a solution
$u$ of Problem \ref{Pf} using Lemma~\ref{l.holder} and Lemma
~\ref{l.limitado}:\\
\textbf{Proof of Lipschitz continuity of $u$:} First, note that it
is sufficient to prove that $u$ is Lipschitz continuous at every
point $y$ on the boundary of $\{u>0\}$, since $\Delta_p u = 0$ on
the open set $\{u>0\}$, and therefore $u$ is
$C_{\textrm{loc}}^{1,\alpha}$ on $\{u>0\}$. Take $y\in
\partial \{u>0\}$ and consider $u_r: B_1 \rightarrow \mathbb{R}$ defined
by:
 $$
 u_r(x) := \frac{1}{r}u(y+rx).
 $$
Since $u$ is a  minimizer of Problem~\ref{Pf} in the ball $B_r(y)$,
we obtain $u_r$ is a local minimizer of the same problem in the ball
$B_1$, with $u_r(0)=u(y)=0$. By Lemma~\ref{l.limitado}, we have that
$\|u_r\|_{L^\infty(B_{1/4})}<C$, where $C$ is a constant depending
only on $p, \gamma$ and $n$. Thus,  by Lemma~\ref{l.holder} we
conclude $u_r \in C^\alpha(B_{1/4})$. Furthermore, for every $r>0$
$$
\|u\|_{C^{\alpha}(B_{1/4})} \leq C=C(n,p, \gamma).
$$
This implies that $\|\nabla u\|_\infty \leq C$ and thus $u$ is
locally $C$-Lipschitz (see Theorem 4.2.3 in \cite{EG92}). This
finish the proof of item (3) of Theorem~\ref{Existence}. $\hfill \square$

At this moment, the proofs of item (4) and (5) are similar to the proofs of Corollary 4.3 and Theorem 4.4 in ~\cite{DP} respectively. The necessary minor modifications are similar to the ones treated in the proof of Lemma \ref{l.limitado} and therefore we will skip the details.
\end{proof}
In the same spirit of \cite{AC}, D. Danielli and A. Petrosyan provided in \cite{DP} a representation theorem for
which Theorem \ref{Existence} puts our minimizers $u_\varepsilon$ under the hypothesis of. The next Theorem will play an important hole in the investigation of fine regularity properties of the free boundary.
\begin{theorem} \label{representation} Let $u_\varepsilon$ be a
minimizer of problem (\ref{Pe}). Then
\begin{enumerate}
\item $\mathcal{H}^{n-1} \big(\mathcal{K} \cap \partial \{ u_\varepsilon
> 0 \}\big ) < \infty$ for every compact set $\mathcal{K} \subset
D^C$.
\item There exists a Borel function $q_\varepsilon$ such that
$$
    \Delta_p u_\varepsilon = q_\varepsilon \mathcal{H}^{n-1}
    \lfloor \partial \{ u_\varepsilon > 0 \},
$$
that is, for any $\psi \in C^\infty_0(D^C)$, there holds
$$
    -\int_{D^C} |\nabla u_\varepsilon|^{p-2}\nabla u_\varepsilon
    \cdot \nabla \psi dx = \int_{ \{ u_\varepsilon > 0 \} } \psi
    q_\varepsilon d \mathcal{H}^{n-1}.
$$
\item For any compact set $\mathcal{K} \subset D^C$, there exist
positive constants $c, ~C$ such that
$$
    c\le q_\varepsilon \le C
$$
and
$$
    cr^{n-1} \le \mathcal{H}^{n-1} \big ( B_r(x) \cap \partial \{
    u_\varepsilon > 0 \} \big ) \le Cr^{n-1},
$$
for every ball $B_x(r) \subset \mathcal{K}$ with $x \in \partial
\{ u_\varepsilon > 0 \}$.
\item For $\mathcal{H}^{n-1}$ almost all points in $\partial \{u_\varepsilon > 0 \}$, an outward normal $\nu = \nu(x)$ is defined and
        furthermore
        $$
            u_\varepsilon(x+y) = q_\varepsilon(x) (y \cdot \nu)^{+} + o(y),
        $$
        where $\frac{o(y)}{|y|} \to 0$ as $ |y| \to 0$. This allows us to define $q_\varepsilon(x) = (u_\varepsilon)_\nu(x)$ at those points.
\item $\mathcal{H}^{n-1} ( \partial \{u_\varepsilon > 0 \} \setminus \partial_\mathrm{red} \{u_\varepsilon > 0 \}) = 0.$
\end{enumerate}
\end{theorem}


\section{A geometric-measure Hadamard's variational formula and the free boundary condition} \label{Hadamard}

In this section we suggest a geometric-measure version of the well
known \linebreak Hadamard's variational formula (see \cite{G}) to
deduce the free boundary condition of Problem \ref{Pe}. Roughly
speaking, given two points in the reduced free boundary, say $x_1$
and $x_2$, the idea is to make an inward perturbation around
$x_1$, an outward perturbation around $x_2$ in such a way that we
do not disturb very much the original volume and then compare the
optimal configuration to the perturbed one in terms of the
functional $ J_\varepsilon$. Here are the details.\par

Let $\rho \colon \mathbb{R} \to \mathbb{R}$ be
a nonnegative $C^\infty$ function supported in $[0,1]$, with, say
$\int \rho(t)dt = 1$. Hereafter, we shall write $u = u_\varepsilon$
and fix two points $x_1$ and $x_2$ in the reduced free boundary
$\partial_\textrm{red} \{ u > 0 \}$. For any $0<r<
\frac{\textrm{dist}(x_1, x_2)}{100}$, and $\lambda>0$, we consider
the vector field
\begin{equation}\label{Perturbation}
    P_r(x) := \left \{
        \begin{array}{cr}
            x + \lambda r \rho \left ( \dfrac{|x-x_1|}{r} \right ) \nu(x_1) & x \in B_r(x_1) \vspace{.2cm}\\
            x - \lambda r \rho \left ( \dfrac{|x-x_2|}{r} \right ) \nu(x_2) & x \in B_r(x_2) \vspace{.2cm} \\
            x & \textrm{ elsewhere.}
        \end{array}
        \right.
\end{equation}

If $\upsilon$ is any vector in $\mathbb{R}^n$, from direct computation, we see that
\begin{equation}\label{Derivative}
DP_r(x) \cdot \upsilon = \upsilon + (-1)^{i+1} \left \{ \lambda \rho'\left ( \dfrac{|x-x_i|}{r} \right ) \dfrac{\langle x-x_i, \upsilon \rangle}{|x-x_i|} \right \} \nu(x_i) ~ \textrm{ in } B_r(x_i)
\end{equation}
Notice that, if $\lambda$ is small enough, $P_r$ is a diffeomorphism that maps $B_r(x_i)$ onto itself. Indeed, if $\lambda \sup\limits_{[0,1]} \rho'(t) < 1$, $P_r$ is a local injective diffeomorphism. Now, if $\lambda \rho(t) \le 1 - t$, for $0\le t \le 1$,
$$
    | P_r(x) - x_i | \le |x - x_i| + \lambda r \rho \left ( \dfrac{|x-x_i|}{r} \right ) \le r,
$$
for any $x \in B_r(x_i)$. Finally, notice that $P_r \Big |_{\partial B_r(x_i)} = Id$, therefore $P_r$ has to be onto.\par
For each $r>0$ small enough, we will consider the $r$-perturbed configuration, $v_r$ implicitly defined by
\begin{equation}\label{vr}
    v_r(P_r(x)) = u(x).
\end{equation}
The idea is to compare our optimal configuration $\{u > 0 \}$ to its perturbation $\{ v_r > 0 \}$ in terms of
the penalized problem \ref{Pe}. An important geometric measure information we shall use is the blow-up limit.
For any $r>0$ small enough and $i=1,2$, consider the blow-up sequence, $ u^i_r \colon B_1(0) \to \mathbb{R}$,
given by
$$
    u^i_r(y) := \dfrac{1}{r} u \big (x_i + ry \big ).
$$
From the blow-up analysis, we know, the set $B_1 \cap \{ u^i_r > 0 \}$ approaches $\{ y \in B_1 ~ \big | ~ \langle y, \nu(x_i) \rangle < 0 \}$, as $r \to 0$. Let us compute the change on the volume of the perturbation. More specifically, make use of the Change of Variables Theorem, we obtain
\begin{equation}\label{Had1}
    \begin{array}{lcl}
        &\dfrac{|\Big ( B_r(x_i) \cap \{ v_r > 0 \} \Big )|}{r^n}&  =  \dfrac{1}{r^n} \displaystyle \int_{ B_r(x_i) \cap \{ v_r > 0 \}} dx \vspace{.2cm}\\
        & = & \displaystyle \int_{ B_1 \cap \{ v_r(x_i + ry) > 0 \}} dy \vspace{.2cm}\\
        & = & \displaystyle \int_{ B_1 \cap \{ u^i_r > 0 \}} \det \left ( DP_r (x_i + r y) \right ) dy \vspace{.2cm}\\
        & \longrightarrow &  \hspace{-.5cm} \displaystyle \int\limits_{ B_1 \cap \{ \langle y, \nu(x_i) \rangle < 0 \} } \hspace{-.5cm} 1 + (-1)^{i+1} \lambda \rho'(|y|) \left \langle \dfrac{y}{|y|}, \nu(x_i) \right \rangle dy,
    \end{array}
\end{equation}
as $r \to 0$. Notice there exists a constant $C(\rho)$ so that, for any unit vector $\nu \in \mathbb{S}^{n-1}$ there holds
\begin{equation}\label{Had2}
C(\rho) \equiv \displaystyle \int\limits_{ B_1 \cap \{ \langle y, \nu \rangle < 0 \} } \hspace{-.5cm} \rho'(|y|) \left \langle \dfrac{y}{|y|}, \nu \right \rangle dy.
\end{equation}
A similar computation shows that
\begin{equation}\label{Had3}
    \dfrac{|\Big ( B_r(x_i) \cap \{ u > 0 \} \Big )|}{r^n} \longrightarrow  \hspace{-.5cm} \displaystyle \int\limits_{ B_1 \cap \{ \langle y, \nu(x_i) \rangle < 0 \} } \hspace{-.5cm} dy,
\end{equation}
as $r \to 0$.
Combining (\ref{Had1}), (\ref{Had2}) and (\ref{Had3}), we conclude
\begin{equation}\label{Had4}
    \dfrac{\Big (|\{ v_r > 0 \}| \Big ) - \Big (|\{ u > 0 \}| \Big ) }{r^n}    \longrightarrow 0,
\end{equation}
as $r \to 0$. From the Lipschitz continuity of the penalization $f_\varepsilon$, we obtain
\begin{equation}\label{Had5}
 f_\varepsilon \Big ( \big ( |\{ v_r > 0 \} |\big ) \Big ) - f_\varepsilon \Big ( | \big ( \{ u > 0 \}| \big ) \Big ) \le \dfrac{1}{\varepsilon} o(r^n).
\end{equation}
Now we shall turn our attention to the differential of the perturbation on the
$p$-Dirichlet integral.
Initially we observe that
\begin{equation}\label{Had6}
    \begin{array}{lcl}
        \dfrac{1}{r^n} \displaystyle \int_{B_r(x_i)} |\nabla u(x)|^p dx & = & \displaystyle \int_{B_1} |\nabla u^i_r(y)|^p dy\vspace{.2cm} \\
        & = & \displaystyle \int_{B_1 \cap \{ u^i_r > 0 \} } |\nabla u^i_r(y)|^p
        dy,
    \end{array}
 \end{equation}
once $\partial_{\textrm{red}} \{ u_r^i >0\}$ is smooth. Now,
applying twice the Change of Variables Theorem, taking into account
that $P_r$ maps $B_r(x_i)$ diffeomorphically onto itself,
\begin{equation}\label{Had7}
    \begin{array}{lcl}
        \dfrac{1}{r^n} \displaystyle\int\limits_{B_r(x_i)} |\nabla v_r(x)|^p dx & = & \dfrac{1}{r^n} \displaystyle \int\limits_{B_r(x_i)} |DP_r(P_r^{-1}(x))^{-1} \cdot \nabla u (P_r^{-1}(x)) |^p  dx \vspace{.2cm} \\
        & = & \dfrac{1}{r^n} \displaystyle \int\limits_{B_r(x_i)} |DP_r(y)^{-1} \cdot \nabla u (y) |^p |\det \big ( DP_r(y) \big ) | dy \vspace{.2cm} \\
        & = & \displaystyle \int\limits_{B_1 \cap \{ u^i_r > 0 \} } |DP_r(x_i + rz)^{-1} \cdot \nabla u^i_r (z) |^p |\det \big ( DP_r(x_i + rz) \big ) | dz. \vspace{.2cm} \\
    \end{array}
\end{equation}
Now, from (\ref{Derivative}), using the fact that for any matrix $A$, with $|A|<1$, we have $\big ( Id + A \big )^{-1} = Id + \sum\limits_{i=1}^\infty (-1)^i A^i$, we have
\begin{equation}\label{Had8}
DP_r(x_i + rz)^{-1} \cdot \nabla u^i_r (z)  = \nabla u^i_r (z) - (1)^{i+1} \lambda \dfrac{\rho'(|z|)}{|z|} \langle z, \nabla u^i_r (z) \rangle \nu(x_i) + o(\lambda).
\end{equation}
On the other hand,
\begin{equation}\label{Had9}
|\det \big ( DP_r(x_i + rz) \big )|  = 1 + (-1)^{i+1}\lambda \dfrac{\rho'(|z|)}{|z|} \langle z,\nu(x_i) \rangle.
\end{equation}
Combining (\ref{Had6}), (\ref{Had7}), (\ref{Had8}) and (\ref{Had9}), we obtain
\begin{equation}\label{Had10}
    \begin{array}{lll}
        & & \dfrac{1}{r^n} \displaystyle \int\limits_{B_r(x_i)} |\nabla v_r(x)|^p - |\nabla u(x)|^p dx  =  (-1)^{i+1}\lambda  \displaystyle \int\limits_{B_1 \cap \{ u^i_r > 0 \} } |\nabla u_r^i(z)|^p\dfrac{\rho'(|z|)}{|z|} \langle z,\nu(x_i) \rangle dz\\
        & + & (-1)^{i}\lambda \displaystyle \int\limits_{B_1 \cap \{ u^i_r > 0 \} } p|\nabla u_r^i(z)|^{p-2}\dfrac{\rho'(|z|)}{|z|}  \langle z, \nabla u_r^i(z) \rangle \langle \nabla u_r^i(z),\nu(x_i) \rangle dz + o(\lambda).
    \end{array}
\end{equation}
Again, from the blow-up analysis (see \cite{DP}), for each $\delta >
0$, we know
$$
    \nabla u_r^i \to q(x_i) \nu(x_i),
$$
uniformly in $B_1 \cap \{ \langle y, \nu(x_i) \rangle < - \delta \}$. Therefore, by $r$-uniform Lipschitz continuity of $u^i_r$, we have
\begin{equation}\label{Had11}
    \nabla u^i_r \to - q(x_i) \nu(x_i) \mathbf{\chi}_{ B_1 \cap \{ \langle y, \nu(x_i) \rangle < 0 \} },
\end{equation}
in $L^p(B_1)$. Letting $r \to 0$ in (\ref{Had10}), we find
\begin{equation}\label{Had12}
\begin{split}
    & \dfrac{1}{r^n} \displaystyle \int_{B_r(x_i)} |\nabla v_r(x)|^p - |\nabla u(x)|^p dx  \longrightarrow \\
    &     (-1)^{i+1}(p-1)\lambda \big (q(x_i)\big )^p \displaystyle \int\limits_{B_1 \cap \{ u^i_r > 0 \} } \dfrac{\rho'(|z|)}{|z|} \langle z,\nu(x_i) \rangle dz + o(\lambda).
\end{split}
\end{equation} Notice that
$$
    \textrm{div}\big ( \rho(|z|) \big ) = \dfrac{\rho'(|z|)}{|z|} \langle z,\nu(x_i) \rangle,
$$
Thus, from Divergence Theorem and the blow-up analysis,
\begin{equation}\label{Had13}
    \displaystyle \int\limits_{B_1 \cap \{ u^i_r > 0 \} } \dfrac{\rho'(|z|)}{|z|} \langle z,\nu(x_i) \rangle dz \rightarrow \int\limits_{B_1 \cap \{ \langle z, \nu(x_i) \rangle =  0 \} } \rho(|z|) d\mathcal{H}^{n-1}(z) = c(\rho).
\end{equation}
Putting (\ref{Had12}) and (\ref{Had13}) together, we obtain
\begin{equation}\label{Had14}
    \displaystyle \int_{D^C} |\nabla v_r(x)|^p - |\nabla u(x)|^p dx = r^n \lambda (p-1)c(\rho)\big ( q(x_1)^p - q(x_2)^p \big ) + r^n
    o(\lambda).
\end{equation}
From the minimality property of $u$, (\ref{Had5}) and (\ref{Had14}),
\begin{equation}\label{Had15}
    0\le J_\varepsilon(v_r) - J_\varepsilon(u) \le r^n \lambda (p-1)c(\rho)\big ( q(x_1)^p - q(x_2)^p \big ) + r^n o(\lambda) + \dfrac{1}{\varepsilon} o(r^n)
\end{equation}
Dividing (\ref{Had15}) by $r^n$ and letting $r \to 0$ we obtain
\begin{equation}\label{Had16}
0 \le \lambda (p-1)c(\rho)\big ( q(x_1)^p - q(x_2)^p \big ) + o(\lambda)
\end{equation}
Now dividing (\ref{Had16}) by $\lambda$, letting $\lambda \to 0$, and afterwards reversing the places of $x_1$ and $x_2$, we finally obtain
\begin{equation}\label{FBC}
    q(x_1) = q(x_2)
\end{equation}
Since $x_1$ and $x_2$ were taking arbitrarily in $\partial_{\textrm{red}} \{ u > 0 \}$, we have proven
\begin{theorem}\label{normal} There exists a positive constant
$\lambda_\varepsilon$ such that
$$
    q_\varepsilon \equiv \lambda_\varepsilon, ~\forall
x \in \partial_{\textrm{red}} \{ u_\varepsilon > 0 \}.
$$
\end{theorem}
It now follows from \cite{DP} that for each $\varepsilon > 0$ fixed, the reduced free boundary $\partial_{\textrm{red}} \{ u_\varepsilon > 0 \}$ is a $C^{1,\alpha}$ smooth surface. Real analyticity of the reduced free boundary is then a consequence of \cite{KNS}. It is worthwhile to point out that, for $n = 2$, a small variant of the main result in \cite{DP1} assures full regularity of the free boundary $ \{ u_\varepsilon > 0 \}$, as long as $p > 2 -\sigma$, for some universal constant $\sigma$. \par

In \cite{C1}, \cite{C2} and \cite{C3}, L. Caffarelli introduced and developed the, by now, well known notion of viscosity solution of a given free boundary problem (for the Laplacian operator). Our final goal of this section is to establish the free boundary condition obtained in Theorem \ref{normal} in the viscosity sense. This will be used as a geometric tool in the remaining sections.

\begin{theorem}[Free boundary condition in the viscosity sense] \label{viscosity} Let $x_0 \in \partial \{u_\varepsilon > 0 \}$ be a free boundary point. Suppose there exists a touching ball $B$, i.e., $\partial B \cap \partial \{ u_\varepsilon > 0 \} = \{x_0\}$, such that either $B \subset \{u_\varepsilon > 0\}$ or $B \subset \{u_\varepsilon = 0 \}$. Then
\begin{equation}\label{viscosity FBC}
    u_\varepsilon(x) = \lambda_\varepsilon \langle x-x_o, \nu \rangle^{+} + o(|x-x_0|),
\end{equation}
where $\lambda_\varepsilon$ is the positive constant provided in Theorem \ref{normal} and $\nu$ is the unit normal vector to $\partial B$, pointing inward to $\{u_\varepsilon > 0 \}$.
\end{theorem}
\begin{proof} Let $B = B_r(\xi)$. We shall first deal with the hypothesis that $B \subset \{u_\varepsilon > 0\}$. For notation convenience, we will omit the center of the ball $B$. With no loss of generality, we can assume $x_0 = 0$ and $\nu = e_n$. Let $\Theta$ be the auxiliary function solving
\begin{equation}
\left \{
    \begin{array}{rcll}
        \textrm{div} \big ( |\nabla u_\varepsilon|^{p-2} D\Theta \big ) &=& 0 & \textrm{ in } \Delta = B_r\setminus \overline{B_{r/2}}\\
        \Theta &=& 1 & \textrm{ on } \partial B_{r/2}\\
        \Theta &=& 0 & \textrm { on } \partial B_{r}
    \end{array}
\right.
\end{equation}

By the nondegeneracy of the gradient $\nabla u_\varepsilon$, see \cite{DP}, $\Theta \in C^{1,\alpha}(\overline{\Delta})$ and since it vanishes on $ \partial B$, we have for some constant $C= C(n,p,\varepsilon, r)$
$$
    \Theta(x) = C x_{n} + o(|x|).
$$
Consider $\tilde{B}$ to be a ball centered at $0$ such that $B \subset\subset \tilde{B}$.  We define
$$
    \theta_{0} = \sup \left\{ m: u_\varepsilon(x)\geq m \Theta (x)  \textrm{ in } \tilde{B}\cap B \right\}.
$$
and any for $k\geq 1$,
$$
    \theta_{k} = \sup \left\{ m: u_\varepsilon (x) \geq m \Theta(x)  \textrm{ in } 2^{-k}\tilde{B} \cap B \right\}.
$$
The sequence $\left\lbrace \theta_{k}\right\rbrace _{k\geq 1}$ is increasing and bounded by a constant $K_\varepsilon$, since $u_\varepsilon$ is Lipschitz continuous. Let $\tilde{\theta} = \sup_{k} \theta_{k}$. If we set
$$
    \theta := C  \tilde{\theta}
$$
we have
$$
    u_\varepsilon (x) \geq \tilde{\theta}\Theta(x) + o(|x|) = \theta x_{N} + o(|x|).
$$
We claim,
\begin{equation} \label{claim viscosity}
    u_\varepsilon (x) = \theta x_{N} + o(|x|).
\end{equation}
Indeed, assume there exists a sequence $(x_{k})_{k\geq 1}$ such that $|x_{k}| = r_{k} \to 0 $ and  for some $\delta_{0} > 0 $,
$$
    u_\varepsilon(x_{k}) - \tilde{\theta}\Theta(x_{k}) > \delta_{0}|x_{k}|.
$$
From the definition of $\theta_{k}$, we can choose a subsequence ${k_{j}}$, where $2^{-6k_{j}} < r_{k_{j}} \leq 2^{-4k_{j}}$ and
\begin{equation}\label{sequence}
    u_\varepsilon(x_{k_{j}}) - \tilde{\theta}_{k_{j}}\Theta(x_{k_{j}}) > \delta_{0}|x_{k_{j}}|.
\end{equation}
Now, the function
$$
    \alpha(x) = u(x) - \theta_{k_{j}}\Theta(x)
$$
satisfies $\textrm{div} \big (|\nabla u_\varepsilon|^{p-2} D\alpha) = 0$ and $\alpha \ge 0$ in $ 2^{-k_{j}}\tilde{B}\cap B$. Since $u_\varepsilon$ is Lipschitz continuous, from (\ref{sequence}), we conclude
$$
    u_\varepsilon(x) - \theta_{k_{j}}h(x) \geq C\delta_{0}|x_{k_{j}}|
$$
in a fixed portion of $\partial B_{r_{k_{j}}}$. From the Poisson representation formula, we have, in, $\frac{1}{2}\tilde{B}_{{r_{k_{j}}}}\cap B$
\begin{equation}
    u_\varepsilon(x) - \theta_{k_{j}}\Theta(x) \geq \overline{C}\delta_{0}|x_{k_{j}}|.
\end{equation}
Now for $j$ large enough, we have
$$
    u_\varepsilon (x) \geq \left( \theta_{k_{j}} + \overline{\delta_{0}}\right) \Theta(x)
$$
which is a contradiction with the definition of $\theta_{k_{j}}$.
We have proven so far, $u_\varepsilon (x) = \theta x_{N} + o(|x|).$ It remains to show, $\theta = \lambda_\varepsilon$. Consider then the blow-up sequence
\begin{equation}\label{blow-up visc}
    u_\varepsilon^k(x) = \dfrac{1}{\rho_k} u_\varepsilon(x_0 + \rho_k x),
\end{equation}
as $\rho_k \to 0$. By Lipschitz continuity and (\ref{claim viscosity}), we know $u_\varepsilon^k$ converges uniform in compact subsets to
$$
    u^\infty_\varepsilon(x) = \theta \langle x, \nu \rangle^{+}.
$$
However, by a small modification of a classical argument, see \cite{AC} or \cite{Fr}, we can show the limit of a blow-up sequence of $u_\varepsilon$ is an absolute minimizer of $J_\varepsilon$ in any ball. Since the free boundary of $ u^\infty_\varepsilon$ is smooth, we obtain from Theorem \ref{normal}, $\theta = \lambda_\varepsilon$, as desired. \par
The case $B_r \subset \{u_\varepsilon = 0\}$, follows by a small modification of the above arguments with the fact that Lemma A1 in \cite{C3} also holds for nonnegative Lipschitz functions $v$ satisfying $\textrm{div}\big( |\nabla u_\varepsilon|^{p-2} Dv \big ) \le 0$.
\end{proof}

\section{Recovering the original problem}
In this section we shall relate a solution to the penalized problem to a (possible) solution to our original problem. Roughly speaking the idea is that the function $f_\varepsilon$ will charge a lot for those configurations that have a volume bigger than 1. We hope if the charge is too big, i.e., if $\varepsilon > 0$ is small enough, optimal configurations of Problem \ref{Pe} will rather prefer to have volume 1 than paying for the penalization. We will follow  the lines of \cite{T}.

\begin{lemma} \label{bound for |u>0|} There exist positive constants $c$ and $C$, independent of $\varepsilon$,
 such that
$$
    c \le | \{ u_\varepsilon > 0 \} | \le 1 + C\varepsilon
$$
\end{lemma}
\begin{proof}
Let $D^\star$ be any smooth domain containing $D$, so that $|D^\star \setminus D| = 1$. From the minimality of
$u_\varepsilon$, we have
\begin{equation}\label{Rec1}
    J_\varepsilon(u_\varepsilon) = \int_{D^C} |\nabla u_\varepsilon(x)|^p + f_\varepsilon \big ( | \{ u_\varepsilon > 0 \} | \big ) \le J_\varepsilon(u^\star) = C,
\end{equation}
where $u^\star$ is the $p$-harmonic function in $D^\star \setminus
D$ taking boundary data equal to $\varphi$ on $\partial D$ and $0$
on $\partial D^\star$. Thus
$$
    \dfrac{1}{\varepsilon} \big ( | \{ u_\varepsilon > 0 \} | - 1 \big ) \le f_\varepsilon \big ( | \{ u_\varepsilon > 0 \} | \big ) \le C.
$$
This proves the estimate from above. Let us turn our attention to the estimate by below. Expression (\ref{Rec1}), together with Poincar\'e inequality, provides
\begin{equation}\label{Rec2}
    \int_{D^C} |\nabla u_\varepsilon(x)|^p + |u_\varepsilon(x)|^p dx \le C,
\end{equation}
for some $C$ independent of $\varepsilon$. Let $D_\delta$ be a tubular neighborhood of $\partial D$. For each  $x_0 \in \partial D$ fixed, let us consider the fiber $F_{x_0} := \{x_0 + t \mu(x_0) ~ \big | ~  0\le t \le \delta \}$ and denote $\Theta_{x_0} := \mathcal{H}^{1}\big ( \{ u_\varepsilon > 0 \} \cap  F_{x_0} \big )$. From Mean Value Inequality, followed by H\"older Inequality and (\ref{Rec2}) we have
\begin{equation}\label{Rec3}
    \begin{array}{lcl}
        \delta \varphi(x_0) &\le& \displaystyle \int_0^\delta u\big ( x_0 + t \mu(x_0) \big ) dt + \displaystyle \int_0^\delta |\nabla u\big ( x_0 + \overline{t} \mu(x_0) \big )|t dt \\
        & \le & C \Theta_{x_0}^{1/q} \delta^{1/q} \left ( 1 + \dfrac{\delta}{q+1} \right ).
    \end{array}
\end{equation}
Now we integrate (\ref{Rec3}) over $\partial D$ and obtain
\begin{equation}\label{Rec4}
    \int_{\partial D} \varphi dS \le C(\delta) | \{ u_\varepsilon > 0 \} \cap D_\delta |^{1/q}.
\end{equation}
Finally, from (\ref{Rec4}), there must exist a constant, independent of $\varepsilon$, so that $ | \{ u_\varepsilon > 0 \} | \ge c$, as claimed.
\end{proof}
\begin{lemma} \label{q < C} There exists a positive constant $C$ independent of $\varepsilon$, so that $\lambda_\varepsilon \le C$, where $\lambda_\varepsilon$ is the constant provided by Theorem \ref{normal}.
\end{lemma}
\begin{proof}   Applying Divergence Theorem to the field $F_1 = u |\nabla u|^{p-2} \nabla u$, we have
\begin{equation}\label{DivThm1}
    \int_{D^C} |\nabla u_\varepsilon|^pdx = \int_{\partial D} \varphi |\nabla u|^{p-2} \partial_\mu u dS,
\end{equation}
where $\mu$ is the outward unit vector in $\partial D$. If we apply Divergence Theorem to the field $F_2 = |\nabla u|^{p-2} \nabla u$, we obtain
\begin{equation}\label{DivThm2}
\lambda_\varepsilon^{p-1} \mathcal{H}^{n-1}(\partial \{ u_\varepsilon > 0 \}) = \int_{\partial D} |\nabla u|^{p-2} \partial_\mu u dS.
\end{equation}
Isoperimetric inequality gives a universal bound by below to $\mathcal{H}^{n-1}(\partial \{ u_\varepsilon > 0 \})$, i.e,  $\mathcal{H}^{n-1}(\partial \{ u_\varepsilon > 0 \}) \ge c$, for some $c$ independent of $\varepsilon$. Combining this with (\ref{DivThm1}) and (\ref{DivThm2}), we obtain
$$
    \lambda_\varepsilon \le C(D, \varphi),
$$
as claimed.
\end{proof}
\begin{lemma}\label{q>c} There exists a universal positive constant $c>0$, such that $ \lambda_\varepsilon \ge c$,
for all $\varepsilon >0$.
\end{lemma}
\begin{proof}
   Let  $z_1 \in D^C$ be such that $u_\varepsilon (z_1) > 0$ for all $\varepsilon>0$. Let us denote by
    $\delta = \textrm{dist}(z_1,\partial D)$. Consider the smooth family of domains $\Upsilon_t := B_{\frac{\delta}{2} +
    t}(z_1) \bigcap D^C$. Let $t_\varepsilon$ denote the first $t$ such that
    $\Upsilon_t$ touches $\partial \{ u_\varepsilon > 0\}$. Let us call $x_0 = \partial \Upsilon_{t_\varepsilon} \bigcap
    \partial \{ u_\varepsilon > 0 \}$. Define $\Psi_\varepsilon$ to be $p$-harmonic function in $\Upsilon_{t_\varepsilon}
    \setminus \Upsilon_0$, with the following boundary values data:
    $$
        \Psi_\varepsilon \big |_{\partial \Upsilon_0} = \min\limits_{\partial D} \varphi  \hspace{.3cm} \textrm{
        and } \hspace{.3cm}  \Psi_\varepsilon \big |_{\partial \Upsilon_{t_\varepsilon}} = 0.
    $$
    By the maximum principle we have $u_\varepsilon \ge \Psi_\varepsilon$ in $\Upsilon_{t_\varepsilon} \setminus \Upsilon_0$. From Hopf's Lemma (see for instance \cite{V}) we also know there exists a constant $c>0$ depending on $\partial D$ and $\inf \varphi$, but independent of $\varepsilon$, such that
\begin{equation}\label{Psi > c}
    \Psi_{-\nu}(x_0) \ge c,
\end{equation}
where $\nu$ denotes the outward unit normal vector of $B_{\frac{\delta}{2} + t_\varepsilon}(z_1)$, at $x_0$. Recall from Theorem \ref{viscosity}, we have the following asymptotic development around $x_0$
\begin{equation} \label{visc}
    \Psi(x) \le u(x) = \lambda_\varepsilon \langle x-x_0, \nu \rangle^{+} + o(|x-x_0|).
\end{equation}
Dividing (\ref{visc}) by $|x-x_0|$, letting $x \to x_0$ and taking into account (\ref{Psi > c}), we finally obtain
$$
    c \le \lambda_\varepsilon,
$$
as desired.
\end{proof}
We are ready to show the main theorem of this section.
\begin{theorem}\label{Pe = P} If $\varepsilon$ is small enough, then any solution of Problem (\ref{Pe}) is a
solution of Problem (\ref{P}).
\end{theorem}

\begin{proof} Let us initially suppose $|\{ u_\varepsilon > 0 \}| > 1$. In the same spirit of Section \ref{Hadamard}, consider a inward perturbation of the set $\{
u_\varepsilon> 0 \}$ with volume change $V$, in such a way that the set of positivity of the new function, $\widetilde{u}_\varepsilon$ is still bigger than $1$. Thus
\begin{equation}\label{vol1}
f_\varepsilon ( | \{ \widetilde{u}_\varepsilon > 0 \}|) -  f_\varepsilon ( | \{ u_\varepsilon > 0 \}|) =
    -\dfrac{1}{\varepsilon} V.
\end{equation}
From (\ref{Had10}) and Lemma \ref{q < C}, we have
\begin{equation} \label{vol2}
    \begin{array}{lll}
        \displaystyle\int_{D^C} | \nabla \widetilde{u}_\varepsilon |^p - | \nabla u_\varepsilon|^p &=& \lambda_\varepsilon^p V + o(V) \\
        &\le & C^p V + o(V).
    \end{array}
\end{equation}
Using the fact that $J_\varepsilon(u_\varepsilon) \le J_\varepsilon(\widetilde{u}_\varepsilon)$, (\ref{vol1}) and (\ref{vol2}), we find
\begin{equation}\label{vol3}
    0 \le C^pV + o(V) - \dfrac{1}{\varepsilon}V.
\end{equation}
Finally if we divide inequality (\ref{vol3}) by $V$ and let $V \to 0$, we obtain
$$
    \varepsilon > \dfrac{1}{C^p}.
$$
If $|\{ u_\varepsilon >0  \}| < 1$, we argue similarly, making an outward perturbation and using Lemma \ref{q>c} to obtain another lower bound for $\varepsilon$. Thus, if $\varepsilon$ is small enough, $|\{u_\varepsilon>0\}|$ automatically adjusts to be equal to $1$.
\end{proof}


\section{Radial Symmetry}
In this section we show a simple symmetry result of Problem (\ref{P}). Indeed, we shall show the best way of insulating a uniformly heated spherical body is by a ball. Recall, when  $\varphi \equiv \textrm{Constant}$, Problem (\ref{P}) is equivalent to our original physical optimization problem.  Here is the theorem:
\begin{theorem} Let $D$ be the unit ball and $\varphi \equiv 1$. Then Problem (\ref{P}) has a unique solution and it is radially symmetric. In particular the free boundary is a sphere.
\end{theorem}
\begin{proof}
Let $u = u_\varepsilon$ be a solution to Problem (\ref{P}), with $D = B_1$ and $\varphi \equiv 1$. Denote $\Omega = \{u > 0 \}$. Let $B_{r_1}$ and $B_{r_2}$ be the biggest ball inside $\Omega \setminus D$ and the smallest ball outside $\Omega$, respectively. Let $y_1 \in \partial B_{r_1} \cap \partial \Omega$ and $y_2 \in \partial B_{r_2} \cap \partial \Omega$. Consider, $h_i$, $i =1,2$ solution of
\begin{equation}\label{pharmonic}
    \left \{
        \begin{array}{rll}
            \Delta_p h_i &= & 0  \textrm{ in } B_{r_i} \setminus B_1 \\
            h_i & = & 1 \textrm{ on } \partial B_1 \\
            h_i & = & 0 \textrm{ on } \partial B_{r_i}. \\
        \end{array}
    \right.
\end{equation}
It is simple to show $h_i$ is radially symmetric. Indeed, $h_i$ is
the unique minimizer of
$$
    E_p(f) := \int_{B_{r_i} \setminus B_1} |\nabla f(x)|^pdx,
$$
among all functions $f \in W^{1,p}$ satisfying the according boundary data. For any orthonormal transformation $\mathcal{O} \in O(n)$, consider
$$
    h_i^\mathcal{O}(x) := h_i(\mathcal{O} x).
$$
Clearly, $h_i^\mathcal{O}$ has the same boundary data as $h_i$ and furthermore,
$$
    E_p(h_i^\mathcal{O}) = \int_{B_{r_i} \setminus B_1} |\mathcal{O}^T \nabla h_i (\mathcal{O} x)|^p dx = E_p(h_i).
$$
Thus $h_i^\mathcal{O} = h_i$. Since $\mathcal{O} \in O(n)$ was
taken arbitrarily, $h_i$ has to be radial. \par In particular, the
inward normal derivative of $h_i$ over $\partial B_{r_i}$ is a
positive constant $\lambda_i$. Since $r_1 \le r_2,$ we have
\begin{equation}\label{ineq1}
    \lambda_1 \ge \lambda_2.
\end{equation}
Now, from the maximum principle,
$$
    h_1 \le u \le h_2.
$$
Hence, from the free boundary condition in the viscosity sense, Theorem \ref{viscosity}, we obtain
\begin{equation}\label{ineq2}
    \lambda_2 = (h_2)_{\nu} (y_2) \le \lambda_\varepsilon \le (h_1)_\nu(y_1) = \lambda_1.
\end{equation}
Combining (\ref{ineq1}) and (\ref{ineq2}), we conclude $\lambda_1 = \lambda_2$, and therefore, $r_1 = r_2$. This implies $\partial \Omega$ has to be a sphere of radius $r_1 = r_2$. \par
We have proven any solution to Problem (\ref{P}), with $D = B_1$ and $\varphi \equiv 1$ is radially symmetric. Uniqueness now follows due to the volume constraint.
\end{proof}


\section*{Acknowledgment}
    The authors would like to thank Professor Irene Gamba (UT-Austin) for having raised the main physical questions that motivate this present work. The second author would like to thank the hospitality of the Universidade Federal de Alagoas, where this work was partially
    developed. Both authors are grateful to Fapeal and Pronex-Dynamical Systems (CNPq) by the financial support.


\bibliographystyle{amsplain, amsalpha}

\end{document}